\documentclass[a4paper,12pt,draft]{amsart}
\usepackage{amsmath, amssymb}
\usepackage[vcentermath,enableskew]{youngtab}
\newcommand{\sgn}{\operatorname{sgn}}
\newcommand{\nodemimus}{\circleddash}
\newcommand{\minimum}{\varnothing}
\newcommand{\partitionof}{\vdash}
\newcommand{\ZZ}{\mathbb{Z}}
\newcommand{\N}{\mathbb{N}}

\newcommand{\der}{\partial}
\newcommand{\Y}{\mathbb{Y}}
\newcommand{\Trees}{\mathbb{T}}
\newcommand{\numof}[1]{|#1|}

\newcommand{\SU}[3][{}]{s^{#2,#3}_{#1}}
\newcommand{\SD}[3][{}]{s_{#2,#3}^{#1}}

\newcommand{\defit}[1]{{\em #1}}

\theoremstyle{plain}
\newtheorem{thm}{Theorem}[section]
\newtheorem{lemma}[thm]{Lemma}
\newtheorem{cor}[thm]{Corollary}
\newtheorem{theorem}[thm]{Theorem}
\newtheorem{proposition}[thm]{Proposition}

\theoremstyle{remark}
\newtheorem{remark}[thm]{Remark}

\newtheorem{example}[thm]{Example}

\theoremstyle{definition}
\newtheorem{definition}[thm]{Definition}

\allowdisplaybreaks[4]

\begin{document}

\title[Pieri's formula]{ Pieri's Formula \\ for Generalized Schur Polynomials}
\author[Numata, y.]{Numata, yasuhide}
\thanks{Department of Mathematics, Hokkaido University, Kita 10, Nishi 8, Kita-Ku, Sapporo,
Hokkaido, 060-0810, Japan.  E-mail: {\tt nu@math.sci.hokudai.ac.jp}}
\address{Affiliation of author: Department of Mathematics, Hokkaido University, Kita 10, Nishi 8, Kita-Ku, Sapporo, Hokkaido, 060-0810, Japan.}
\begin{abstract}
Young's lattice, the lattice of all Young diagrams,
 has the Robinson-Schensted-Knuth correspondence,
the correspondence between
 certain matrices and pairs of 
semi-standard  Young tableaux with the same shape.
Fomin introduced generalized Schur operators  
 to generalize the  Robinson-Schensted-Knuth correspondence.
In this sense,
generalized Schur operators are generalizations 
of semi-standard Young tableaux.
We define a generalization of Schur polynomials as 
 expansion coefficients  of  generalized Schur operators.
We show that
the commutating relation of generalized Schur operators
 implies Pieri's formula to 
 generalized Schur polynomials.
\end{abstract}

\maketitle

\section{Introduction}
Young's lattice is a prototypical example of a differential poset
which was first defined by Stanley \cite{StaD, StaV}.
The Robinson correspondence is
a correspondence between permutations and
pairs of standard tableaux whose shapes are the same Young diagram.
This correspondence was generalized for differential posets or
dual graphs (generalizations of differential posets \cite{FomD}) 
by Fomin \cite{FomRSK, FomS}. (See also \cite{Rob}.)

Young's lattice also has  
The Robinson-Schensted-Knuth correspondence, 
the correspondence between
 certain matrices and pairs of semi-standard tableaux.
Fomin \cite{FomK}
introduced operators called generalized Schur operators, and
generalized 
the Robinson-Schensted-Knuth correspondence
for generalized Schur operators.
We define a generalization of Schur polynomials 
as  expansion coefficients  of  generalized Schur operators.

A complete symmetric polynomial is a Schur polynomial associated with 
a Young diagram consisting of only one row.
Schur polynomials satisfy Pieri's formula,
the formula describing  
the product of a complete symmetric polynomial and a Schur polynomial as 
a sum of Schur polynomials$:$
\begin{gather*}
h_i(t_1,\ldots,t_n)s_\lambda(t_1,\ldots,t_n)=\sum_\mu
 s_\mu(t_1,\ldots,t_n),
\end{gather*}
where the sum is over all $\mu$'s that are obtained from $\lambda$
by adding $i$ boxes, with no two in the same column,
$h_i$ is the $i$-th complete symmetric polynomial, and 
$s_\lambda$ is the Schur polynomial associated with $\lambda$.

In this paper, we generalize 
 Pieri's formula to generalized Schur polynomials.

\begin{remark}
 Lam introduced a generalization of the Boson-Fermion
 correspondence \cite{lam}.
In the paper, he also showed Pieri's and Cauchy's formulae for 
some families of symmetric functions in the 
context of Heisenberg algebras. 
Some important families of symmetric functions,
 e.g., Schur functions, Hall-Littlewood
 polynomials, Macdonald polynomials and so on, 
are examples of them.
He proved Pieri's formula using essentially the same method as the one
 in this paper. 
Since the assumptions of generalized Schur operators are less than
those of  Heisenberg algebras,
 our polynomials  are more general than his; 
e.g., some of our polynomials are 
not symmetric. 
An example of generalized Schur operators which provides non-symmetric
 polynomials is in Section \ref{extree}. 
See also Remark \ref{lamsresult}
for the relation between \cite{lam} and this paper.
\end{remark}

This paper is organized as follows:
In Section $\ref{GSdefsec}$, we recall generalized Schur operators, and
define generalized Schur polynomials. We also define
a generalization of complete symmetric polynomials, called weighted
complete symmetric polynomials, in Section $\ref{WCdefsec}$.
In Section $\ref{mainsec}$, we show Pieri's formula for these
polynomials (Theorem $\ref{genepieri}$).
We also see that Theorem $\ref{genepieri}$ becomes simple for special
parameters, and that weighted complete symmetric polynomials are written
as linear combinations of generalized
Schur polynomials in a special case. 
Other examples are shown in  Section $\ref{EXsec}$.

\section{Definition}

We introduce two types of polynomials in this section.
One is a generalization of Schur polynomials.
The other  is a generalization of complete symmetric polynomials.

\subsection{Generalized Schur Polynomials}\label{GSdefsec}

First we recall the generalized Schur  operators defined by Fomin \cite{FomK}.
We define a generalization of Schur polynomials
as expansion coefficients of generalized Schur operators.

Let $K$ be a field of characteristic zero
that contains all formal power series in variables
$t, t', t_1,t_2,\ldots$ 
Let 
$V_i$ be finite-dimensional $K$-vector spaces for all $i \in \ZZ$. 
Fix a basis $Y_i$ of each $V_i$
so that $V_i=KY_i$.
Let $Y=\bigcup_i Y_i$,  $V=\bigoplus_i V_i$ and $\widehat{V}=\prod_i V_i$,  
i.e., $V$ is the vector space consisting of all finite linear
combinations of elements of $Y$ and
$\widehat{V}$ is the vector space consisting of all linear combinations
of elements of $Y$.
The \defit{rank function} on $V$  mapping $v \in V_i$ to $i$ 
is denoted by $\rho$. 
We say that $Y$  \defit{has a minimum} $\minimum$
if $Y_i=\emptyset$ for $i<0$ and $Y_0=\{\minimum\}$.

For a sequence $\{A_i\}$ and a formal variable $x$, 
we write $A(x)$ for the generating function $\sum_{i\geq 0} A_i x^i$.

\begin{definition}
Let $D_i$ and $U_i$ be linear maps on $V$
for nonnegative integers $i\in \N$.
We call $D(t_1)\cdots D(t_n)$ and $U(t_n)\cdots U(t_1)$
\defit{generalized Schur operators with} $\{a_m\}$
if the following conditions are satisfied$:$
\begin{itemize}
\item $\{a_m\}$ is a sequence of  $K$.
\item $U_i$  satisfies
$U_i(V_j) \subset V_{j+i}$
for all $j$.
\item $D_i$ satisfies
$D_i(V_j) \subset V_{j-i}$
for all $j$.
\item The equation 
$D(t')U(t)=a(t t')U(t)D(t')$
holds.
\end{itemize}
\end{definition}

In general, $D(t_1)\cdots D(t_n)$ and $U(t_n)\cdots U(t_1)$ are
not linear operators on $V$  but  linear operators from $V$ to $\widehat{V}$.

Let  $\langle\phantom{x},\phantom{x}\rangle$ 
be the natural pairing,
i.e.,
the bilinear form 
on $\widehat{V} \times V$
such that 
$\langle\sum_{\lambda\in Y}a_\lambda \lambda,\sum_{\mu \in Y}b_\mu \mu\rangle =\sum_{\lambda \in Y}a_{\lambda} b_{\lambda}$.
For generalized Schur operators $D(t_1)\cdots D(t_n)$ and $U(t_n)\cdots U(t_1)$,
$U_i^{\ast}$ and $D_i^{\ast}$ denote
the maps obtained from the adjoints of $U_i$ and $D_i$ 
with respect to  $\langle\phantom{x},\phantom{x}\rangle$
by restricting to $V$, respectively.
For all $i$, 
$U_i^{\ast}$ and $D_i^{\ast}$ 
are linear maps on $V$
satisfying $U_i^{\ast}(V_j)\subset V_{j-i}$
and $D_i^{\ast}(V_{j})\subset V_{j+i}$.
It follows by definition that
\begin{align*}
  \langle v, U_{i} w \rangle & = \langle w, U_{i}^{\ast} v \rangle,
&  \langle v, D_{i} w \rangle & = \langle w, D_{i}^{\ast} v \rangle
\end{align*}
for $v$, $w\in V$.
We write $U^{\ast}(t)$ and $D^{\ast}(t)$ for 
$\sum U^{\ast}_i t^i$ and
$\sum D^{\ast}_i t^i$.
It follows by definition that
\begin{align*}
  \langle  U(t) \mu ,\lambda\rangle &= \langle  U^{\ast}(t) \lambda , \mu\rangle,
&  \langle  D(t) \mu ,\lambda\rangle &= \langle  D^{\ast}(t) \lambda , \mu\rangle
\end{align*}
for $\lambda$, $\mu\in Y$.
The equation  $D(t')U(t)=a(t t')U(t)D(t')$ implies the equation
$U^{\ast}(t')D^{\ast}(t)=a(t t')D^{\ast}(t)U^{\ast}(t')$.
Hence  $U^{\ast}(t_1)\cdots U^{\ast}(t_n)$ and 
$D^{\ast}(t_n)\cdots D^{\ast}(t_1)$ are 
generalized Schur operators with $\{a_m\}$ when 
 $D(t_1)\cdots D(t_n)$ and $U(t_n)\cdots U(t_1)$
are.

\begin{definition}
Let $D(t_1)\cdots D(t_n)$ and $U(t_n)\cdots U(t_1)$
be generalized Schur operators with $\{a_m\}$.
For $v \in V$ and $\mu \in Y$,
$\SD[D]{v}{\mu}(t_1,\ldots,t_n)$
 and $\SU[U]{\mu}{v}(t_1,\ldots,t_n)$
are respectively defined by
\begin{align*}
\SD[D]{v}{\mu}(t_1,\ldots,t_n)&=\langle D(t_1)\cdots D(t_n) v,\mu\rangle,\\  
\SU[U]{\mu}{v}(t_1,\ldots,t_n)&=\langle U(t_n)\cdots U(t_1) v,\mu\rangle.
\end{align*}
We call these polynomials $\SD[D]{v}{\mu}(t_1,\ldots,t_n)$
 and $\SU[U]{\mu}{v}(t_1,\ldots,t_n)$ 
\defit{generalized Schur polynomials}.
\end{definition}

\begin{remark}
Generalized Schur polynomials $\SD[D]{v}{\mu}(t_1,\ldots,t_n)$
are symmetric 
in the case when $D(t)D(t')=D(t')D(t)$,
but not symmetric in general.
Similarly,
generalized Schur polynomials $\SU[U]{\mu}{v}(t_1,\ldots,t_n)$
are symmetric 
if $U(t)U(t')=U(t')U(t)$.

If $U_0$ (resp. $D_0$) is the identity map on $V$, 
generalized Schur polynomials $\SD[D]{v}{\mu}(t_1,\ldots,t_n)$
(resp. $\SU[U]{\mu}{v}(t_1,\ldots,t_n)$) 
are quasi-symmetric.
In \cite{BMSW},
Bergeron, Mykytiuk, Sottile and van Willigenburg 
considered graded representations of the algebra of
 noncommutative symmetric functions on the $\ZZ$-free module 
whose basis is a graded poset,
and
gave a Hopf-morphism
from a Hopf algebra generated by intervals of the poset to 
the Hopf algebra of quasi-symmetric functions.
\end{remark}

\begin{example}\label{protex}
Our prototypical example  is
Young's lattice $\Y$ that consists of all Young diagrams.
Let $Y$ be Young's lattice $\Y$,
 $V$ the $K$-vector space $K\Y$ whose basis is $\Y$,
and $\rho$ the ordinary rank function
mapping a Young diagram $\lambda$ to the number of boxes in $\lambda$.
Young's lattice $\Y$ has a minimum $\minimum$, 
 the Young diagram with no boxes.
We call a skew Young diagram $\mu/\lambda$ a horizontal strip
if $\mu/\lambda$ has no two boxes in the same column.
Define $U_i$ 
by $U_i(\mu)=\sum_{\lambda}\lambda$, 
where
the sum is over all $\lambda$'s that are
obtained from $\mu$ by adding a horizontal strip consisting of
$i$ boxes;
and define $D_i$ by
$D_i(\lambda)=\sum_{\mu}\mu$, 
where
the sum is over all $\mu$'s that are
obtained from $\lambda$ by removing a horizontal strip consisting of $i$ boxes.
For example, 
\begin{align*}
\Yboxdim{6pt}\yng(3,1,1)\ &\mathop{\longmapsto}^{D_2}\ 
\Yboxdim{6pt}\yng(2,1) + \Yboxdim{6pt}\yng(1,1,1) \\
\Yboxdim{6pt}\yng(2,1)  \ &\mathop{\longmapsto}^{U_2}\ 
\Yboxdim{6pt}\yng(2,2,1) + \Yboxdim{6pt}\yng(3,1,1) + 
\Yboxdim{6pt}\yng(3,2) + \Yboxdim{6pt}\yng(4,1).
\end{align*}
(See also Figure $\ref{figofy}$, the graph of $D_1$ ($U_1$) and $D_2$ ($U_2$).)  
\begin{figure}
\begin{picture}(220,135)(-110,-32)\scriptsize
\put(0,-25){\makebox(0,0)[c]{$\minimum$}}
\put(0,0){\makebox(0,0)[c]{$(1)$}}
\put(25,25){\makebox(0,0)[c]{$(2)$}}
\put(-25,25){\makebox(0,0)[c]{$(11)$}}
\put(50,50){\makebox(0,0)[c]{$(3)$}}
\put(0,50){\makebox(0,0)[c]{$(21)$}}
\put(-50,50){\makebox(0,0)[c]{$(111)$}}
\put(75,75){\makebox(0,0)[c]{$(4)$}}
\put(25,75){\makebox(0,0)[c]{$(31)$}}
\put(0,75){\makebox(0,0)[c]{$(22)$}}
\put(-25,75){\makebox(0,0)[c]{$(211)$}}
\put(-75,75){\makebox(0,0)[c]{$(1111)$}}
\put(100,100){\makebox(0,0)[c]{$(5)$}}
\put(50,100){\makebox(0,0)[c]{$(41)$}}
\put(25,100){\makebox(0,0)[c]{$(32)$}}
\put(0,100){\makebox(0,0)[c]{$(311)$}}
\put(-25,100){\makebox(0,0)[c]{$(221)$}}
\put(-50,100){\makebox(0,0)[c]{$(211)$}}
\put(-100,100){\makebox(0,0)[c]{$(11111)$}}
\put(0,-18){\line(0,1){10}}
\put(-7,7){\line(-1,1){10}}
\put( 7,7){\line(1,1){10}}
\put(32,32){\line(1,1){10}}
\put(18,32){\line(-1,1){10}}
\put(-18,32){\line(1,1){10}}
\put(-32,32){\line(-1,1){10}}
\put(57,57){\line(1,1){10}}
\put(43,57){\line(-1,1){10}}
\put(7,57){\line(1,1){10}}
\put(0,57){\line(0,1){10}}
\put(-57,57){\line(-1,1){10}}
\put(-43,57){\line(1,1){10}}
\put(-7,57){\line(-1,1){10}}
\put(82,82){\line(1,1){10}}
\put(68,82){\line(-1,1){10}}
\put(32,82){\line(1,1){10}}
\put(25,82){\line(0,1){10}}
\put(18,82){\line(-1,1){10}}
\put(7,82){\line(1,1){10}}
\put(-7,82){\line(-1,1){10}}
\put(-82,82){\line(-1,1){10}}
\put(-68,82){\line(1,1){10}}
\put(-32,82){\line(-1,1){10}}
\put(-25,82){\line(0,1){10}}
\put(-18,82){\line(1,1){10}}
\qbezier[30](7,-25)(25,0)(25,18)
\qbezier[30](7,0)(50,0)(50,43)
\qbezier[20](0,7)(0,25)(0,43)
\qbezier[30](32,25)(75,25)(75,68)
\qbezier[20](25,32)(25,50)(25,68)
\qbezier[30](-15,25)(23,25)(23,68)
\qbezier[20](-25,32)(-25,50)(-25,68)
\qbezier[30](57,50)(100,50)(100,93)
\qbezier[20](50,57)(50,75)(50,93)
\qbezier[30](-38,50)(-7,75)(-2,93)
\qbezier[20](-50,57)(-50,75)(-50,93)
\qbezier[30](4,57)(13,75)(21,93)
\qbezier[30](2,57)(15,75)(2,93)
\qbezier[30](-4,57)(-13,75)(-21,93)
\put(95,-10){\makebox(0,0)[r]{--- $D_1$ ($U_1$)}}
\put(95,-25){\makebox(0,0)[r]{$\cdots$ $D_2$ ($U_2$)}}
\end{picture}
\caption{Young's lattice}\label{figofy}
\end{figure}

In this case, $D(t_1)\cdots D(t_n)$ and $U(t_n)\cdots U(t_1)$
are generalized Schur operators
with $\{a_m=1\}$.
Both
$\SD[D]{\lambda}{\mu}(t_1,\ldots,t_n)$
 and $\SU[U]{\lambda}{\mu}(t_1,\ldots,t_n)$
are equal to the skew Schur polynomial $s_{\lambda/\mu}(t_1,\ldots,t_n)$
for $\lambda, \mu \in \Y$.
For example, since
\begin{align*}
D(t_2) \Yboxdim{6pt}\yng(2,1) = & \Yboxdim{6pt}\yng(2,1) + t_2\Yboxdim{6pt}\yng(1,1) + t_2\Yboxdim{6pt}\yng(2) + t_2^2\Yboxdim{6pt}\yng(1) \\ 
D(t_1)D(t_2) \Yboxdim{6pt}\yng(2,1)= &\Yboxdim{6pt}\yng(2,1) + t_1\Yboxdim{6pt}\yng(1,1) + t_1\Yboxdim{6pt}\yng(2) + t_1^2\Yboxdim{6pt}\yng(1) \\
&+ t_2(\Yboxdim{6pt}\yng(1,1) + t_1\Yboxdim{6pt}\yng(1)) + t_2(\Yboxdim{6pt}\yng(2) + t_1\Yboxdim{6pt}\yng(1) +t_1^2 \minimum) \\
&+ t_2^2(\Yboxdim{6pt}\yng(1)+t_1 \minimum) ,
\end{align*}
$\SD[D]{(2,1)}{\minimum}(t_1,t_2)=s_{(2,1)}(t_1,t_2)=t_1^2 t_2+t_1 t_2^2$.
\end{example}

\begin{example}\label{secex}
Our second example is the polynomial ring $K [x]$ with a
variable $x$. Let $V$ be $K[x]$ and 
 $\rho$ 
 the ordinary rank function  mapping
a monomial $ax^n$ to its degree $n$.
In this case, $\dim V_i=1$ for all $i \geq 0$ and $\dim V_i=0$ for $i< 0$.
Hence its basis $Y$ is identified with $\N$
and has a minimum $c_0$,
a nonzero constant.
Define $D_i$ and $U_i$ by $\frac{\der^i}{i!}$
and $\frac{x^i}{i!}$,
where $\der$ is the partial differential operator in $x$.
Then $D(t)$ and $U(t)$ are $\exp(t\der)$ and $\exp(t x)$.
Since $D(t)$ and $U(t)$ satisfy $D(t)U(t')=\exp(t t')U(t')D(t)$,
 $D(t_1)\cdots D(t_n)$ and $U(t_n)\cdots U(t_1)$ are 
generalized Schur operators with $\{a_m=\frac{1}{m!}\}$.
In general, for  differential posets,  we can construct
generalized Schur operators in a similar manner.

Since $\der$ and $x$ commute with $t$,
the following equations hold${:}$
\begin{gather*}
D(t_1)\cdots D(t_n)=\exp(\der t_1)\cdots\exp(\der t_n)
=\exp(\der(t_1+\cdots+t_n)),\\
U(t_n)\cdots U(t_1)=\exp(x t_n)\cdots\exp(x t_1)
=\exp(x(t_1+\cdots+t_n)).
\end{gather*}
It follows from direct calculations that
\begin{align*}
\exp(\der(t_1+\cdots+t_n))c_{i}x^{i}
&=\sum_{j=0}^{i}\frac{(t_1+\cdots+t_n)^j}{j!}\frac{i!}{(i-j)!}c_ix^{i-j}\\
&=\sum_{j=0}^{i}\frac{{i!}(t_1+\cdots+t_n)^j c_i}{{(i-j)!j!}c_{i-j}}c_{i-j}x^{i-j},\\
\exp(x(t_1+\cdots+t_n))c_ix^i
&=\sum_j\frac{(t_1+\cdots+t_n)^j x^j}{j!}c_ix^i\\
&=\sum_j\frac{(t_1+\cdots+t_n)^j c_i}{j!c_{i+j}}c_{i+j}x^{i+j}.
\end{align*}
Hence it follows that
\begin{align*}
\SD[D]{c_{i+j}x^{i+j}}{c_{i}x^i}(t_1,\ldots,t_n)&=\frac{(i+j)!}{i!j!}\frac{c_{i+j}}{c_i}(t_1+\ldots+t_n)^j\\
\SU[U]{c_{i+j}x^{i+j}}{c_{i}x^i}(t_1,\ldots,t_n)&=\frac{1}{j!}\frac{c_{i}}{c_{i+j}}(t_1+\ldots+t_n)^j,
\end{align*}
if we take $\{c_ix^i\}$ as the basis $Y$.

If $c_i=1$ for all $i$, then 
$\SD[D]{x^{i+j}}{x^i}(t_1,\ldots,t_n)=\frac{(i+j)!}{i!j!}(t_1+\ldots+t_n)^j$,
and $\SU[U]{x^{i+j}}{x^i}(t_1,\ldots,t_n)=\frac{1}{j!}{(t_1+\ldots+t_n)^j}$.
\end{example}

\begin{lemma}\label{dualitylemma}
Generalized Schur polynomials satisfy the following equations{\rm :} 
\begin{gather*}
\SD[D]{\lambda}{\mu}(t_1,\ldots,t_n)
=\SU[D^{\ast}]{\lambda}{\mu}(t_1,\ldots,t_n),\\
\SU[U]{\lambda}{\mu}(t_1,\ldots,t_n)
=\SD[U^{\ast}]{\lambda}{\mu}(t_1,\ldots,t_n)
\end{gather*}
for $\lambda, \mu\in Y$.
Generalized Schur polynomials also satisfy the following equations$:$
\begin{gather*}
\SD[D]{v}{\mu}(t_1,\ldots,t_n)
=\sum_{\nu\in Y} \langle v,\nu\rangle
\SU[D^{\ast}]{\nu}{\mu}(t_1,\ldots,t_n),\\
\SU[D^{\ast}]{\mu}{v}(t_1,\ldots,t_n)=
\sum_{\nu\in Y} \langle v,\nu \rangle \SD[D]{\mu}{\nu}(t_1,\ldots,t_n)
,\\
\SD[U^{\ast}]{v}{\mu}(t_1,\ldots,t_n)=
\sum_{\nu\in Y}\langle v,\nu\rangle\SU[U]{\nu}{\mu}(t_1,\ldots,t_n),\\
\SU[U]{\mu}{v}(t_1,\ldots,t_n)
=\sum_{\nu\in Y}\langle v,\nu\rangle \SD[U^{\ast}]{\mu}{\nu}(t_1,\ldots,t_n)
\end{gather*}
for $\mu \in Y$, $v \in V$.
\end{lemma}
\begin{proof}
It follows by definition that
\begin{align*}
\SD[D]{\lambda}{\mu}(t_1,\ldots,t_n)
&=\langle D(t_1)\cdots D(t_n) \lambda,\mu \rangle\\
&=\langle D^{\ast}(t_n)\cdots D^{\ast}(t_1) \mu , \lambda\rangle
=\SU[D^{\ast}]{\lambda}{\mu}(t_1,\ldots,t_n).
\end{align*}
Similarly, we have 
$\SU[U]{\lambda}{\mu}(t_1,\ldots,t_n)=\SD[U^{\ast}]{\lambda}{\mu}(t_1,\ldots,t_n)$.
The other formulae follow from 
$v=\sum_{\nu\in Y}\langle\nu, v\rangle\nu$ for $v\in V$.
\end{proof}

\begin{remark}
Rewriting the generalized Cauchy identity \cite[1.4. Corollary]{FomK}
with our notation,
we obtain a Cauchy identity for generalized Schur polynomials:
\begin{align*}
\sum_{\nu\in Y}
\SD[D]{\nu}{\mu}&(t_1,\ldots,t_n)
\SU[U]{\nu}{v}(t'_1,\ldots,t'_n)\\
&=
\prod_{i,j} a(t_i t'_j)
\sum_{\kappa \in Y}
\SU[U]{\mu}{\kappa}(t'_1,\ldots,t'_n)
\SD[D]{v}{\kappa}(t_1,\ldots,t_n)
\end{align*}
for $v\in V$, $\mu\in Y$.
\end{remark}

\begin{remark}\label{lamsresult}
In this remark, we construct operators $B_l$ from generalized Schur
 operators  $D(t_1)\cdots D(t_n)$ and $U(t_n)\cdots U(t_1)$.
These operators $B_l$ are closely related to the results of Lam \cite{lam}.
Furthermore we can construct  other generalized
Schur operators $D(t_1)\cdots D(t_n)$ and $U'(t_n)\cdots U'(t_1)$ 
from $B_l$.

Let   $D(t_1)\cdots D(t_n)$ and $U(t_n)\cdots U(t_1)$ be
generalized Schur operators with  $\{a_m\}$.
For a partition $\lambda \partitionof l$,
we define $z_\lambda$  by
 $z_\lambda=1^{m_1(\lambda)}m_1(\lambda)!\cdot 2^{m_2(\lambda)}m_2(\lambda)!\cdots$, where
$m_i(\lambda)=\numof{\{j|\lambda_j=i\}}$.
Let $U_0=D_0=I$, where $I$ is the identity map.
For positive integers $l$, we inductively define   
$b_l$, $B_l$ and  $B_{-l}$ by
\begin{align*}
b_l=&a_l-\sum_{\lambda} \frac{b_\lambda}{z_\lambda}, \\
B_l=&D_l-\sum_{\lambda} \frac{ B_\lambda}{z_\lambda}, \\
B_{-l}=&U_{l}-\sum_{\lambda}  \frac{ B_{-\lambda}}{z_\lambda},
\end{align*}
where
$b_\lambda=b_{\lambda_1}\cdot b_{\lambda_2}\cdots$,
$B_\lambda=B_{\lambda_1}\cdot B_{\lambda_2}\cdots$,
$B_{-\lambda}=B_{-\lambda_1}\cdot B_{-\lambda_2}\cdots$ and
 the sums are over all partitions $\lambda$ of  $l$ such that
$\lambda_1 < l$. Let $b_l\neq 0$ for any $l$.
It follows from direct calculations that
\begin{align*}
[B_{l},B_{-l}] &= l\cdot b_l \cdot I,\\
[B_{l},B_{-k}] &= 0 
\end{align*}
for positive integers $l \neq k$.
If $U_i$ and $D_i$ respectively commute with $U_j$ and $D_j$ for all $i,j$,
then $\{B_{l},B_{-l}| l\in \ZZ_{>0}\}$ generates the Heisenberg algebra.
In this case, we can apply the results of Lam \cite{lam}.
See also Remark \ref{remcomp} for the relation between his complete symmetric
 polynomials $h_i[b_m](t_1,\ldots,t_n)$ and our weighted complete
 symmetric polynomials $h_i^{\{a_m\}}(t_1,\ldots,t_n)$.

For a partition $\lambda \partitionof l$,
let $\sgn(\lambda)$ denote $(-1)^{\sum_i (\lambda_i - 1)}$, where 
the sum is over all $i$'s such that $\lambda_i > 0$.
Although  $U_i$ and $D_i$ do not commute with $U_j$ and $D_j$, we can define dual generalized Schur operators 
 $D(t_1)\cdots D(t_n)$ and $U'(t_n)\cdots U'(t_1)$ with
 $\{a'_m\}$ by 
\begin{align*}
a'_l=&\sum_{\lambda} \frac{\sgn(\lambda) b_\lambda}{z_\lambda}, \\
U'_{-l}=&\sum_{\lambda}  \frac{\sgn(\lambda) B_{-\lambda}}{z_\lambda},
\end{align*}
where  the sums are over all partitions $\lambda$ of  $l$.
In this case, it follows from direct calculations that $a(t) \cdot a'(-t) = 1$.
\end{remark}

\subsection{Weighted Complete Symmetric Polynomials}\label{WCdefsec}
Next we introduce a generalization of complete symmetric polynomials.
We define weighted symmetric polynomials inductively.

\begin{definition}
Let $\{a_m\}$ be a sequence of elements of $K$.
We define the $i$-th weighted complete symmetric polynomial
$h^{\{a_m\}}_{i}(t_1,\ldots,t_n)$
to be the coefficient of $t^{i}$ in $a(t_1 t)\cdots a(t_n t)$.
\end{definition}

By definition, 
for each $i$,
the $i$-th weighted complete symmetric polynomial
$h^{\{a_m\}}_{i}(t_1,\ldots,t_n)$
is a homogeneous symmetric polynomial of degree $i$.

\begin{remark}
For a  sequence $\{a_m\}$ of elements of $K$, 
the $i$-th weighted complete symmetric polynomial
$h_{i}^{\{a_m\}}(t_1,\ldots,t_n)$ coincides with the polynomial
defined by
\begin{gather}
h_{i}^{\{a_m\}}(t_1,\ldots,t_n)=
\begin{cases}
a_it^i_1 & \text{(for $n=1$)},\\
\displaystyle\sum_{j=0}^{i}h_{j}^{\{a_m\}}(t_1,\ldots,t_{n-1})h_{i-j}^{\{a_m\}}(t_n)\label{eq1}
& \text{(for $n>1$)}.
\end{cases}
\end{gather}
\end{remark}

\begin{example}\label{comex}
When $a_m$ equals $1$ for each $m$,
 $a(t)=\sum_it^i=\frac{1}{1- t}$.
In this case, $h_{j}^{\{1,1,\ldots\}}(t_1,\ldots,t_n)$ equals 
the complete symmetric polynomial $h_j(t_1,\ldots,t_n)$.
\end{example}

\begin{example}
When $a_m$ equals $\frac{1}{m!}$ for each $m$,
$\sum_{j}h_{j}^{\{\frac{1}{m!}\}}(t) = \exp(t) =a(t)$ and
$h_{j}^{\{\frac{1}{m!}\}}(t_1,\ldots,t_n)=\frac{1}{j!}(t_1+\cdots+t_n)^j$.
\end{example}

\begin{remark}\label{remcomp}
In this remark, we compare
 the complete symmetric polynomials $h_i [b_m](t_1,\ldots,t_n)$ of Lam \cite{lam}
and our weighted complete symmetric polynomials $h_i^{\{a_m\}}(t_1,\ldots,t_n)$.
Let $\{b_m\}$ be a sequence of elements of $K$.
The polynomials $h_i [b_m](t_1,\ldots,t_n)$ of  Lam are defined by
\begin{gather*}
h_i[b_m](t_1,\ldots,t_n)=
\sum_{\lambda \partitionof i} \frac{b_\lambda p_\lambda (t_1,\ldots,t_n)}{z_\lambda},
\end{gather*}
where $b_\lambda=b_{\lambda_1}\cdot b_{\lambda_2}\cdots$,
$p_\lambda(t_1,\ldots,t_n)=p_{\lambda_1}(t_1,\ldots,t_n)\cdot p_{\lambda_2}(t_1,\ldots,t_n)\cdots$
 and
$p_i(t_1,\ldots,t_n)=t_1^i+\cdots+t_n^i$.
These polynomials satisfy the equation
\begin{gather*}
h_i[b_m](t_1,\ldots,t_n)=
\sum_{j=0}^{i}h_{j}[b_m](t_1,\ldots,t_{n-1})h_{i-j}[b_m](t_n).
\end{gather*}
Let $a_i=\sum_{\lambda \partitionof i} \frac{b_\lambda}{z_\lambda} $.
Then it follows $h_i[b_m](t_1)= a_i t^i$.
Hence 
\begin{gather*}
h_i[b_m](t_1,\ldots,t_n)=h_i^{\{a_m\}}(t_1,\ldots,t_n).
\end{gather*}
\end{remark}

\section{Main Results}\label{mainsec}

In this section, we show 
some properties of generalized Schur
polynomials and weighted complete
symmetric polynomials.

Throughout this section, let 
 $D(t_1)\cdots D(t_n)$ and  $U(t_n)\cdots U(t_1)$ be generalized
 Schur operators with $\{a_m\}$.

\subsection{Main Theorem}
In Proposition \ref{comutarel},
 we describe the commuting relation of $U_i$ and $D(t_1)\cdots D(t_n)$, 
proved in Section $\ref{proofsec}$.
This relation implies Pieri's formula for our polynomials (Theorem $\ref{genepieri}$),
the main result in this paper.
It also follows from this relation that the weighted complete
symmetric polynomials are written as linear combinations of generalized
Schur polynomials when $Y$ has a minimum (Proposition \ref{meanofcompsym}).
 
First  we describe the commuting relation of $U_i$ and $D(t_1)\cdots
D(t_n)$.
We  prove it in Section $\ref{proofsec}$.
\begin{proposition}\label{comutarel} \label{comutarel2}
The equations 
\begin{align}
D(t_1)\cdots D(t_n)U_i
&=\sum_{j=0}^{i}h_{i-j}^{\{a_m\}}(t_1,\ldots,t_n)U_j D(t_1)\cdots D(t_n)
,\label{eqfirst}\\
D_i U(t_n)\cdots U(t_1)
&=\sum_{j=0}^{i}
h_{i-j}^{\{a_m\}}(t_1,\ldots,t_n)U(t_n)\cdots U(t_1)D_j
,\label{eq37}\\
U^{\ast}_i D^{\ast}(t_n)\cdots D^{\ast}(t_1)
&=\sum_{j=0}^{i}
h_{i-j}^{\{a_m\}}(t_1,\ldots,t_n)D^{\ast}(t_n)\cdots D^{\ast}(t_1)U^{\ast}_j
,\label{eq36}\\
U^{\ast}(t_1)\cdots U^{\ast}(t_n)D^{\ast}_i
&=\sum_{j=0}^{i}
h_{i-j}^{\{a_m\}}(t_1,\ldots,t_n)
D^{\ast}_j U^{\ast}(t_1)\cdots U^{\ast}(t_n)\label{eq38}.
\end{align}
hold for all $i$. 
\end{proposition}

These equations imply the following main theorem.

\begin{theorem}[Pieri's formula]\label{genepieri}
For each $\mu\in Y_k$ and each $v \in V$,
generalized Schur polynomials satisfy
\begin{gather*}
\SD[D]{U_i v}{\mu}
(t_1,\ldots,t_n)
=\sum_{j=0}^{i}
h_{i-j}^{\{a_m\}}(t_1,\ldots,t_n)\sum_{\nu\in Y_{k-j}}
\langle U_j\nu , \mu \rangle
\SD[D]{v}{\nu}(t_1,\ldots,t_n).
\end{gather*}
\end{theorem}
\begin{proof}
It follows from Proposition \ref{comutarel} 
that 
\begin{align*}
\langle D(t_1)\cdots D(t_n)U_i v,\mu\rangle 
&=
\langle \sum_{j=0}^{i}
h_{i-j}^{\{a_m\}}(t_1,\ldots,t_n)U_j D(t_1)\cdots D(t_n) v
,\mu\rangle\\
&=
 \sum_{j=0}^{i}
h_{i-j}^{\{a_m\}}(t_1,\ldots,t_n)\langle U_j D(t_1)\cdots D(t_n) v
,\mu\rangle
\end{align*}
 for $v \in V$ and $\mu \in Y$.
This says 
\begin{align*}
\SD[D]{U_i v}{\mu}&(t_1,\ldots,t_n)\\
&=\sum_{j=0}^{i}
h_{i-j}^{\{a_m\}}(t_1,\ldots,t_n)\sum_{\nu\in Y_{k-j}}
\langle U_j\nu , \mu \rangle
\SD[D]{v}{\nu}(t_1,\ldots,t_n).
\end{align*}
\end{proof}

This formula becomes simple in the case when
$v \in Y$.

\begin{cor}
For each $\lambda, \mu\in Y$,
generalized Schur polynomials satisfy
\begin{gather*}
\SD[D]{U_i \lambda}{\mu}(t_1,\ldots,t_n)
=\sum_{j=0}^{i}
h_{i-j}^{\{a_m\}}(t_1,\ldots,t_n) \cdot 
\SU[D^{\ast}]{\lambda}{U^{\ast}_j \mu}(t_1,\ldots,t_n).
\end{gather*}
\end{cor}
\begin{proof}
It follows from Theorem $\ref{genepieri}$
that
\begin{gather*}
\SD[D]{U_i \lambda}{\mu}(t_1,\ldots,t_n)
=\sum_{j=0}^{i}
h_{i-j}^{\{a_m\}}(t_1,\ldots,t_n)
\sum_{\nu\in Y}
\langle U_j\nu , \mu \rangle
\SD[D]{\lambda}{\nu}(t_1,\ldots,t_n).
\end{gather*}
Lemma $\ref{dualitylemma}$ implies 
\begin{align*}
\sum_{\nu\in Y}
\langle U_j\nu , \mu \rangle
\SD[D]{\lambda}{\nu}(t_1,\ldots,t_n)
&=\sum_{\nu\in Y}
\langle \nu , U_j^{\ast}\mu \rangle
\SD[D]{\lambda}{\nu}(t_1,\ldots,t_n)\\
&=\SU[D^{\ast}]{\lambda}{U^{\ast}_j \mu}(t_1,\ldots,t_n).
\end{align*}
Hence
\begin{gather*}
\SD[D]{U_i \lambda}{\mu}(t_1,\ldots,t_n)
=\sum_{j=0}^{i}
h_{i-j}^{\{a_m\}}(t_1,\ldots,t_n) \cdot 
\SU[D^{\ast}]{\lambda}{U^{\ast}_j \mu}(t_1,\ldots,t_n).
\end{gather*}
\end{proof}

If $Y$ has a minimum $\minimum$, Theorem $\ref{genepieri}$ implies
the following corollary.

\begin{cor}\label{piericor}
Let $Y$ have a minimum $\minimum$.
For all $v \in V$,
the following equations hold$:$
\begin{align*}
\SD[D]{U_i v}{\minimum}(t_1,\ldots,t_n)&
=u_0 \cdot h_{i}^{\{a_m\}}(t_1,\ldots,t_n)
 \cdot \SD[D]{v}{\minimum}(t_1,\ldots,t_n),
\end{align*}
where
$u_0$ is the element of $K$ that satisfies $U_0\minimum=u_0\minimum$.
\end{cor}

In the case when $Y$ has a minimum $\minimum$, 
 weighted complete symmetric polynomials are written as
linear combinations of generalized Schur polynomials.

\begin{proposition}\label{meanofcompsym}
Let $Y$ have a minimum $\minimum$.
The following
equations hold for all $i\geq 0{:}$
\begin{gather*}
\SD[D]{U_i\minimum}{\minimum}(t_1,\ldots,t_n)=d_0^n u_0  \cdot h_{i}^{\{a_m\}}(t_1,\ldots,t_n),
\end{gather*}
where $d_0$, $u_0$ are the elements of $K$ that satisfy $D_0 \minimum = d_0 \minimum$
and $U_0 \minimum = u_0 \minimum$.
\end{proposition}

\begin{proof}
By definition, 
$\SD[D]{\minimum}{\minimum}(t_1,\ldots,t_n)$
is $d_0^n$.
Hence it follows from Corollary \ref{piericor} that
\begin{gather*}
\SD[D]{U_i\minimum}{\minimum}(t_1,\ldots,t_n) = 
u_0 h_{i}^{\{a_m\}}(t_1,\ldots,t_n) d_0^n .
\end{gather*}
 \end{proof}

\begin{example}\label{protexpieri}
In the prototypical example $\Y$ (Example \ref{protex}),
for $\lambda\in \Y$,
$U_i\lambda$ is the sum of all
Young diagrams obtained from $\lambda$
by adding a horizontal strip consisting of $i$ boxes.
Hence $\SD[D]{U_i \lambda}{\minimum}(t_1,\ldots,t_n)$ equals
$\sum_{\nu} s_\nu$, where 
the sum is over all $\nu$'s that
are obtained from $\lambda$
by adding  a horizontal strip consisting of $i$ boxes.
On the other hand, $u_0$ is $1$, and $h_{i}^{\{1,1,1,\ldots\}}(t_1,\ldots,t_n)$
is the $i$-th complete symmetric polynomial $h_i(t_1,\ldots,t_n)$
(Example $\ref{comex}$).
Thus Corollary $\ref{piericor}$ 
is nothing but  the classical Pieri's formula.
Theorem $\ref{genepieri}$ is Pieri's formula for 
skew Schur polynomials; for a skew Young diagram $\lambda / \mu$
and $i\in\N$,
\begin{gather*}
\sum_{\kappa}s_{ \kappa / \mu}(t_1,\ldots,t_n)
=\sum_{j=0}^{i}
\sum_{\nu} h_{i-j}(t_1,\ldots,t_n)
s_{\lambda/\nu}(t_1,\ldots,t_n),
\end{gather*}
where 
the first sum is over all $\kappa$'s that are obtained from $\lambda$
by adding  a horizontal strip consisting of $i$ boxes;
the last sum is over all $\nu$'s that are obtained from $\mu$
by removing  a horizontal strip consisting of $j$ boxes.

In this example, 
Proposition \ref{meanofcompsym} says that the Schur polynomial $s_{(i)}$ corresponding to 
Young diagram with only one row equals the complete symmetric
polynomial $h_i$.
\end{example}

\begin{example}
In the second example $\N$ (Example \ref{secex}),
Proposition \ref{meanofcompsym} says that 
the constant term of 
$\exp(\der(t_1+\cdots+t_n))\cdot\frac{x^i}{i!}$
equals $\frac{(t_1+\cdots+t_n)^i}{i!}$.
\end{example}

\subsection{Some Variations of Pieri's Formula}
In this section,
we show some variations of Pieri's formula for generalized Schur
polynomials,
i.e., we show Pieri's formula not only for $\SD[D]{\lambda}{\mu}(t_1,\ldots,t_n)$
but also 
for $\SU[U]{\lambda}{\mu}(t_1,\ldots,t_n)$,
 $\SU[D^\ast]{\lambda}{\mu}(t_1,\ldots,t_n)$
and
$\SD[U^\ast]{\lambda}{\mu}(t_1,\ldots,t_n)$.

\begin{theorem}[Pieri's formula]\label{genepieri2}
For each $\mu\in Y_k$ and each $v \in V$,
generalized Schur polynomials satisfy the following equations$:$
\begin{align*}
\sum_{\kappa\in Y}\langle D_i \kappa,\mu\rangle
&\SU[U]{\kappa}{v}(t_1,\ldots,t_n)\\
&=\sum_{j=0}^{i} h_{i-j}^{\{a_m\}}(t_1,\ldots,t_n)
\SU[U]{\mu}{D_j v}(t_1,\ldots,t_n),\\
\SD[U^{\ast}]{D^{\ast}_i v}{\mu}
(t_1,\ldots,&t_n)\\
&=\sum_{j=0}^{i}
h_{i-j}^{\{a_m\}}(t_1,\ldots,t_n)\sum_{\nu\in Y_{k-j}}
\langle D^{\ast}_j\nu , \mu \rangle
\SD[U^{\ast}]{v}{\nu}(t_1,\ldots,t_n),\\
\sum_{\kappa\in Y}\langle U^{\ast}_i \kappa,\mu\rangle
&\SU[D^{\ast}]{\kappa}{v}(t_1,\ldots,t_n)\\
&=\sum_{j=0}^{i} h_{i-j}^{\{a_m\}}(t_1,\ldots,t_n)
\SU[D^{\ast}]{\mu}{U^{\ast}_j v}(t_1,\ldots,t_n).
\end{align*}
\end{theorem}

\begin{proof}
Applying 
Theorem $\ref{genepieri}$ 
to $U^{\ast}(t_1)\cdots U^{\ast}(t_n)$
and $D^{\ast}(t_n)\cdots D^{\ast}(t_1)$,
we obtain
\begin{gather*}
\SD[U^{\ast}]{D^{\ast}_i v}{\mu}(t_1,\ldots,t_n)
=\sum_{j=0}^{i}
h_{i-j}^{\{a_m\}}(t_1,\ldots,t_n)\sum_{\nu\in Y_{k-j}}
\langle D^{\ast}_j\nu , \mu \rangle
\SD[U^{\ast}]{v}{\nu}(t_1,\ldots,t_n).
\end{gather*}

It follows from Proposition \ref{comutarel2} 
that 
\begin{gather*}
\langle  D_i U(t_n)\cdots U(t_1) v, \mu \rangle
 =
\langle 
\sum_{j=0}^{i} h_{i-j}^{\{a_m\}}(t_1,\ldots,t_n)U(t_n)\cdots U(t_1) D_j
 v 
,\mu \rangle
\end{gather*}
 for $v \in V$ and $\mu \in Y$.
This equation says
\begin{gather*}
\sum_{\kappa\in Y}\langle D_i \kappa,\mu\rangle
\SU[U]{\kappa}{v}(t_1,\ldots,t_n)
=\sum_{j=0}^{i} h_{i-j}^{\{a_m\}}(t_1,\ldots,t_n)
\SU[U]{\mu}{D_j v}(t_1,\ldots,t_n).
\end{gather*}
For generalized Schur operators 
$U^{\ast}(t_1)\cdots U^{\ast}(t_n)$ and
$D^{\ast}(t_n)\cdots D^{\ast}(t_1)$,
this equation  says
\begin{gather*}
\sum_{\kappa\in Y}\langle U^{\ast}_i \kappa,\mu\rangle
\SU[D^{\ast}]{\kappa}{v}(t_1,\ldots,t_n)
=\sum_{j=0}^{i} h_{i-j}^{\{a_m\}}(t_1,\ldots,t_n)
\SU[D^{\ast}]{\mu}{U^{\ast}_j v}(t_1,\ldots,t_n).
\end{gather*}

\end{proof}

\begin{cor}\label{piericor2}
For all $v \in V$,
the following equations hold$:$
\begin{align*}
\SD[U^{\ast}]{D^{\ast}_i v}{\minimum}(t_1,\ldots,t_n)
&= d_0 \cdot h_{i}^{\{a_m\}}(t_1,\ldots,t_n)
 \cdot \SD[U^{\ast}]{ v}{\minimum}(t_1,\ldots,t_n),
\end{align*}
where
$d_0$ is the element of $K$ that satisfies
$D_0\minimum=d_0\minimum$.
\end{cor}
\begin{proof}
We obtain 
this proposition from Theorem $\ref{piericor}$ 
by applying to generalized Schur operators $U^{\ast}(t_1,\ldots,t_n)$
and $D^{\ast}(t_1,\ldots,t_n).$
\end{proof}

\begin{proposition}\label{meanofcompsym2}
Let $Y$ have a minimum $\minimum$. Then 
\begin{gather*}
\SD[U^{\ast}]{D^{\ast}_i\minimum}{\minimum}(t_1,\ldots,t_n)
=u_0^nd_0 \cdot h_{i}^{\{a_m\}}(t_1,\ldots,t_n),
\end{gather*}
where $u_0$ and $d_0$ are the elements of $K$ 
that satisfy $D_0 \minimum = d_0 \minimum$
and
 $U_0 \minimum = u_0 \minimum$.
\end{proposition}
\begin{proof}
We obtain this proposition 
by applying  Theorem $\ref{meanofcompsym}$ to generalized Schur
 operators 
$U^{\ast}(t_1)\cdots U^{\ast}(t_n)$
and $D^{\ast}(t_n)\cdots D^{\ast}(t_1).$
\end{proof}

\subsection{Proof of Proposition \ref{comutarel}}\label{proofsec}
In this section, we prove Proposition \ref{comutarel}.

First, we prove the equation  $(\ref{eqfirst})$.
The other equations  follow from  the equation  $(\ref{eqfirst})$.

\begin{proof}
Since $D(t_1)\cdots D(t_n)$ and $U(t_n)\cdots U(t_1)$
are generalized Schur operators 
 with $\{a_m\}$,
the equations $D(t)U_i =\sum_{j=0}^{i}a_j t^j U_{i-j} D(t)$
 hold for all integers $i$.
Hence $D(t_1)\cdots D(t_n)U_i$
is written as a $K$-linear combination of
$U_{j}D(t_1)\cdots D(t_n)$.
We write $H_{i,j}(t_1,\ldots,t_n)$
for the coefficient of $U_{j}D(t_1)\cdots D(t_n)$
in $D(t_1)\cdots D(t_n)U_i$.
 
It follows from the equation $D(t)U_i =\sum_{j=0}^{i}a_j t^j U_{i-j} D(t)$ 
that 
\begin{gather}
H_{i,i-j}(t_1) = a_j t_1^j \label{kome}
\end{gather}
 for $0\leq j \leq i$.

We apply the relation ($\ref{kome}$)
to $D(t_n)$ and $U_i$ to have
\begin{gather*}
D(t_1)\cdots D(t_{n-1})D(t_n)U_i
=\sum_{j=0}^{i}a_{i-j}t_n^{i-j}D(t_1)\cdots D(t_{n-1})U_j D(t_n).
\end{gather*}
Since $D(t_1)\cdots D(t_{n-1}) U_i=\sum_j H_{i,j}(t_1,\ldots,t_{n-1})
U_j D(t_1)\cdots D(t_{n-1})$ by the definition of
$H_{i,j}(t_1,\ldots,t_{n-1})$, 
we have the equation
\begin{align*}
D(t_1)\cdots &D(t_{n-1})D(t_n)U_i\\
=&\sum_{k=0}^{i}\sum_{j=k}^{i}a_{i-j}t_n^{i-j}H_{j,k}(t_1,\ldots,t_{n-1})U_{k}D(t_1)\cdots D(t_n).
\end{align*}

Since
$D(t_1)\cdots D(t_n)U_i$
equals $\sum_{k=0}^{i} H_{i,k}(t_1,\ldots,t_n)U_k D(t_1)\cdots D(t_n)$
by definition, the equation
\begin{gather*}
\sum_{k=0}^{i}\sum_{j=k}^{i}a_{i-j}t_n^{i-j}H_{j,k}(t_1,\ldots,t_{n-1})U_{k}D(t_1)\cdots D(t_n)\\
=\sum_{k=0}^{i} H_{i,k}(t_1,\ldots,t_n)U_k D(t_1)\cdots D(t_n)
\end{gather*}
holds.
Hence the equation
\begin{gather}
\sum_{j=k}^{i}a_{i-j}t_n^{i-j}H_{j,k}(t_1,\ldots,t_{n-1})=
H_{i,k}(t_1,\ldots,t_n) \label{staa}
\end{gather} holds.

We claim that $H_{i+k,k}(t_1,\ldots,t_n)$ does not depend on $k$.
It follows from this relation $(\ref{staa})$ that
\begin{align*}
 H_{k+l,k}(t_1,\ldots,t_n)=&\sum_{j=k}^{k+l}a_{k+l-j}
t_n^{k+l-j}H_{j,k}(t_1,\ldots,t_{n-1})\\
=&\sum_{j'=0}^{l}a_{k+l-(j'+k)}t_n^{k+l-(j'+k)}
H_{j'+k,k}(t_1,\ldots,t_{n-1})\\
=&\sum_{j'=0}^{l}a_{l-j'}t_n^{l-j'}
H_{j'+k,k}(t_1,\ldots,t_{n-1}).
\end{align*}
Since the monomials $a_{l-j'}t_n^{l-j'}$ do not depend on $k$,
the equations  
\begin{gather*}
H_{(i-k)+k,k}(t_1,\ldots,t_n)=H_{(i-k)+k',k'}(t_1,\ldots,t_n)
\end{gather*}
hold
if the equations $H_{k+j,k}(t_1,\ldots,t_{n-1}) = H_{k'+j,k'}(t_1,\ldots,t_{n-1})$ hold
 for all $k$, $k'$ and $j\leq i-k$. 
In fact, 
since  $H_{i+k,k}(t_1)$ equals $a_it_1^i$,
$H_{i+k,k}(t_1)$ does not depend on $k$.
Hence it follows inductively that
 $H_{i+k,k}(t_1,\ldots,t_n)$ does not 
depend on $k$, either.
Hence we may write $\tilde{H}_{i-j}(t_1,\ldots,t_n)$ for $H_{i,j}(t_1,\ldots,t_n)$.

It follows from the
equations $(\ref{kome})$ 
 and $(\ref{staa})$ that 
\begin{gather*}
\begin{cases}
\tilde{H}_{i}(t_1)=a_i t_1^i &\text{(for $n=1$)},\\
\tilde{H}_{i}(t_1,\ldots,t_n)=
\sum_{k=0}^i\tilde{H}_{i-k}(t_1,\ldots,t_{n-1})\tilde{H}_{k}(t_n)
& \text{(for $n>1$)}.
\end{cases}
\end{gather*}
Since
$\tilde{H}_{i}(t_1,\ldots,t_n)$ equals
 the $i$-th weighted complete symmetric polynomial
 $h^{\{a_m\}}_{i}(t_1,\ldots,t_n)$,
 we have the equation  $(\ref{eqfirst})$.

We obtain
the equation $(\ref{eq36})$ from the equation $(\ref{eqfirst})$
by applying $\ast$.

Since  $D(t_1)\cdots D(t_n)$ and $U(t_n)\cdots U(t_1)$
are generalized Schur operators with $\{a_m\}$,
 $U^{\ast}(t_1)\cdots U^{\ast}(t_n)$ and $D^{\ast}(t_n)\cdots D^{\ast}(t_1)$ are also generalized Schur operators with $\{a_m\}$.
Applying the equation $(\ref{eq36})$
Proposition $\ref{comutarel}$ to
$U^{\ast}(t_1)\cdots U^{\ast}(t_n)$ and $D^{\ast}(t_n)\cdots D^{\ast}(t_1)$,
we obtain the equation  $(\ref{eq37})$  and $(\ref{eq38})$, respectively.

Hence Proposition $\ref{comutarel}$ follows.
\end{proof}

\section{More Examples}\label{EXsec}
In this section, we consider some examples of generalized Schur operators.

\subsection{Shifted Shapes}
This example is the same as \cite[Example 2.1]{FomK}.
Let $Y$ be the set of shifted shapes, i.e., 
\begin{gather*}
Y=\left\{
\{(i,j)\in \N^2| i\leq j \leq \lambda_i + i\}
\big| \lambda=(\lambda_1>\lambda_2>\cdots), \lambda_i \in \N\right\}
.
\end{gather*}
For $\lambda\subset \nu \in Y$,
let $cc_0(\lambda\setminus \nu)$ denote the number of connected
components of $\lambda\setminus\nu$ that do not intersect 
with the main diagonal,
and $cc(\lambda\setminus \nu)$  the number of connected
components of $\lambda\setminus\nu$.
For example, let 
$\lambda=(7,5,3,2)$ and $\mu=(5,4,2)$. In this case, $\lambda\setminus\nu$
is the set of boxes $\circ$ and $\bullet$ in 
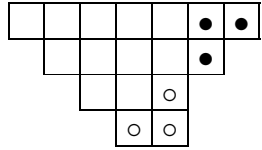
\begin{gather*}
\young(\ \ \ \ \ \bullet \bullet ,:\ \ \ \ \bullet ,::\ \ \circ ,:::\circ \circ ).
\end{gather*}
Since the component of the boxes $\circ$ intersects with the main diagonal 
at $(4,4)$, 
$cc_0(\lambda\setminus \nu)=1$ and $cc(\lambda\setminus \nu)=2$.

For $\lambda\in Y$,  $D_i$
 are defined by
\begin{gather*}
D_i\lambda=\sum_\nu 2^{cc_0(\lambda\setminus\nu)} \nu,
\end{gather*}
where
the sum is over all $\nu$'s 
that are obtained from $\lambda$ by removing $i$ boxes,
with no two box in the same diagonal.

For $\lambda\in Y$, $U_i$
 are defined  by
\begin{gather*}
U_i\lambda=\sum_\mu 2^{cc(\mu\setminus\lambda)} \mu,
\end{gather*}
where
the sum is over all $\mu$'s  
that are obtained from $\lambda$ by
adding $i$-boxes,
with no two box in the same diagonal.
 
In this case, since $D(t)$ and $U(t)$ satisfy
\begin{gather*}
D(t')U(t)=\frac{1+t t'}{1-t t'}U(t)D(t'),
\end{gather*}
$D(t_1)\cdots D(t_n)$ and $U(t_n)\cdots U(t_1)$ are
generalized Schur operators with $\{1,2,2,2,\ldots\}$.
(See \cite{FomK}.)
In this case, for $\lambda$, $\mu \in Y$,
generalized Schur polynomials $\SD[D]{\lambda}{\mu}$
and $\SU[U]{\lambda}{\mu}$ are respectively
the shifted skew Schur polynomials
$Q_{\lambda/\mu}(t_1,\ldots,t_n)$
and
$P_{\lambda/\mu}(t_1,\ldots,t_n)$.

In this case, Proposition  $\ref{meanofcompsym}$ reads as
\begin{gather*}
h^{\{1,2,2,2,\ldots\}}_i (t_1,\ldots,t_n)=
\begin{cases}
 2 Q_{(i)}(t_1,\ldots,t_n) & i>0\\
Q_{\minimum}(t_1,\ldots,t_n) & i=0
\end{cases}
.
\end{gather*}
It also follows from Proposition $\ref{meanofcompsym2}$ that
\begin{gather*}
h^{\{1,2,2,2,\ldots\}}_i (t_1,\ldots,t_n)=
P_{(i)}(t_1,\ldots,t_n).
\end{gather*}

Furthermore, Corollary $\ref{piericor}$
reads as 
\begin{gather*}
\sum_\mu 2^{cc(\mu \setminus \lambda)}Q_\mu(t_1,\ldots,t_n)
=h^{\{1,2,2,2,\ldots\}}_i Q_\lambda(t_1,\ldots,t_n),
\end{gather*}
where
the sum is over all $\mu$'s
that are obtained from 
$\lambda$ by adding $i$ boxes,
with no two in the same diagonal.

\subsection{Young's Lattice: Dual Identities}
This example is the same as \cite[Example 2.4]{FomK}.
Let $Y$ be Young's lattice $\Y$,
and $D_i$ the same ones in the prototypical example, 
(i.e., $D_i \lambda =\sum_\mu \mu$,  where the sum
is over all $\mu$'s that are obtained from $\lambda$
by removing $i$ boxes, with no two in the same column.)
For $\lambda\in Y$, $U'_i$ are defined by
$U'_i\lambda=\sum_\mu \mu$,
where the sum
is over all $\mu$'s that are obtained from $\lambda$
by adding $i$ boxes, with no two in the same row.
(In other words, 
$D_i$ removes  horizontal strips,
while $U'_i$ adds  vertical strips.)

In this case, since $D(t)$ and $U'(t)$ satisfy
\begin{gather*}
D(t)U'(t')=({1+t t'})U'(t')D(t),
\end{gather*}
$D(t_1)\cdots D(t_n)$ and $U'(t_n)\cdots U'(t_1)$ are
generalized Schur operators with $\{1,1,0,0,0,\ldots\}$.
(See \cite{FomK}.)
In this case, for $\lambda$, $\mu \in Y$,
generalized Schur polynomials $\SU[U']{\lambda}{\mu}$
equal $s_{\lambda'/\mu'}(t_1,\ldots,t_n)$,
where 
$\lambda'$ and $\mu'$ are the transposes of $\lambda$ and $\mu$,
and $s_{\lambda'/\mu'}(t_1,\ldots,t_n)$ are skew Schur polynomials.

In the prototypical example (Example \ref{protexpieri}),
 Corollary $\ref{piericor}$ is the classical Pieri's formula, 
the formula describing the product of a complete symmetric polynomial
 and a Schur polynomial. 
In this example,
 Corollary $\ref{piericor}$ is the dual Pieri's formula, 
the formula describing the product of a elementary symmetric polynomial
 and a Schur polynomial.

In this case, Corollary $\ref{meanofcompsym}$  reads as
\begin{gather*}
h^{\{1,1,0,0,0,\ldots\}}_i (t_1,\ldots,t_n)=
 s_{(1^i)}(t_1,\ldots,t_n) = e_i(t_1,\ldots,t_n),
\end{gather*}
where 
$e_i(t_1,\ldots,t_n)$ denotes the $i$-th elementally symmetric polynomials.

Furthermore, Corollary $\ref{piericor}$
 reads as 
\begin{gather*}
\sum_\mu s_\mu(t_1,\ldots,t_n)=
e_i(t_1,\ldots,t_n) s_\lambda(t_1,\ldots,t_n),
\end{gather*}
where
the sum is over all $\mu$'s 
that are obtained from $\lambda$ by adding 
a vertical strip consisting of $i$ boxes.

 For a skew Young diagram $\lambda / \mu$
and $i\in\N$, Theorem $\ref{genepieri}$  reads as
\begin{gather*}
\sum_{\kappa}s_{ \kappa / \mu}(t_1,\ldots,t_n)
=\sum_{j=0}^{i}
\sum_{\nu} e_{i-j}(t_1,\ldots,t_n)
s_{\lambda/\nu}(t_1,\ldots,t_n),
\end{gather*}
where 
the first sum is over all $\kappa$'s that are obtained from $\lambda$
by adding a vertical strip consisting $i$ boxes;
the last sum is over all $\nu$'s that are obtained from $\mu$
by removing a vertical strip consisting $j$ boxes.

\subsection{Planar Binary Trees}\label{extree}
This example is the same as \cite{tree}.
Let $F$ be the monoid of words generated by the alphabet $\{1,2\}$ and
$0$ denote the word of length $0$.
We give  $F$ the structure of a poset by $v \leq v w$ for $v,w \in F$.
We call an ideal of the poset $F$ a \defit{planar binary tree} or shortly a \defit{tree}.
An element of a tree is called a \defit{node} of the tree.
We write $\Trees$
for the set of trees and $\Trees_i$ for the set of trees with $i$ nodes.
We respectively call nodes $v2$ and $v1$ \defit{right} and \defit{left
children} of $v$. 
A node without a child is called a \defit{leaf}.
For $T\in \Trees$ and $v \in F$,
we define $T_v$ to be $\{w\in T| v\leq w\}$.

First we define up operators.
We respectively call $T'$ a tree obtained from $T$ 
by adding some nodes \defit{right-strictly} and \defit{left-strictly}
if $T \subset T'$ and  
 each $w \in T'\setminus T$ has no right children and no left children.
We define linear operators $U_i$ and $U'_i$ on $K \Trees$ by
\begin{gather*}
U_i T = \sum_{T'}  T',\\
U'_i T = \sum_{T''}  T'',
\end{gather*}
where the first sum is over all $T'$'s that are 
obtained from $T$ by adding $i$
 nodes right-strictly, and 
 the second sum is over all $T''$'s that are
obtained from $T$ by adding $i$
 nodes left-strictly.
For example,
\begin{align*}
U_2 \{ 0 \} =& \{0,1,11\} +\{0,1,2\} +\{0,2,21\},\\
U'_2 \{ 0 \}=& \{0,2,22\} +\{0,1,2\} +\{0,1,12\}.
\end{align*}

Next we define down operators.
For $T\in \Trees$, let $r_T$ be $\{ w \in T | w2\not\in T.$ 
If $w=v1w'$ then $v2 \not\in T$. $\}$, i.e.,
the set of nodes which
have no child on its right
and which belong between $0$ and the rightmost leaf of $T$.
The set $r_T$ is a chain.
Let $r_T = \{  w_{T,1} < w_{T,2} < \cdots \}$.
We define linear operators $D_i$ on $K \Trees$ by
\begin{gather*}
D_i T = 
\begin{cases}
(\cdots((T \nodemimus w_{T,i}) \nodemimus w_{T,i-1})\cdots)
 \nodemimus w_{T,1} & i \leq \numof{r_T}\\
0& i >   \numof{r_T}
\end{cases}
\end{gather*}
for $T\in \Trees$, where
\begin{gather*}
T \nodemimus w = (T \setminus T_w) \cup \{w v|w1v \in T_w\}
\end{gather*}
for $w \in T$ such that $w2 \not\in T$.
Roughly speaking, $D_i T $ 
is the tree obtained from $T$ by 
evacuating the $i$ topmost nodes 
without a child on its right and belonging between $0$ and the rightmost leaf of $T$.
For example, let $T$ be $\{0,1,11,12,121\}$. Since $w_{T,1}=0$,
$w_{T,2}=12$ and
\begin{gather*}
\{0,1,11,12,121\} \mathop{\longrightarrow}^{\nodemimus 12} \{0,1,11,12\} \mathop{\longrightarrow}^{\nodemimus 0} \{0,1,2\},
\end{gather*}
 we have $D_2 T= \{0,1,2\}$.

These operators $D(t)$, $U(t')$ and  $U'(t')$ satisfy the following
equations:
\begin{align*}
D(t)U(t')&=\frac{1}{1-t t'}U(t')D(t),\\
D(t)U'(t')&=(1+t t')U'(t')D(t).
\end{align*}
(See  \cite{tree} for a proof of the equations.)
Hence the generalized Schur polynomials for these operators
satisfy the same Pieri's formula as in the case of the classical Young's
lattice and its dual construction.

In this case,  generalized Schur polynomials are not symmetric in general.
For example, since
\begin{align*}
&D(t_1)D(t_2) \{ 0, 1 , 12\}\\
&=D(t_1)( \{ 0, 1 , 12\} + t_2 \{0,2\} + t_2^2 \{0\})\\
&=( \{ 0, 1 , 12\} + t_1 \{0,2\} + t_1^2 \{0\})
+t_2(\{0,2\}+ t_1 \{0\}) 
+t_2^2(\{0\} +t_1\emptyset),
\end{align*}
$\SD[D]{\{0,1,12\}}{\emptyset}(t_1,t_2)=t_1 t_2^2$ is not symmetric.

We define three kinds of labeling on trees to give generalized Schur
polynomials 
$\SU[U]{T}{\emptyset}(t_1,\ldots,t_n)$,
$\SU[U']{T}{\emptyset}(t_1,\ldots,t_n)$
and $\SD[D]{T}{\emptyset}(t_1,\ldots,t_n)$
presentations as generating functions of them.

\begin{definition}
Let $T$ be a tree and $m$ a positive integer.
We call a map $\varphi: T \to \{1,\ldots , m\}$
a \defit{right-strictly-increasing labeling} if 
\begin{itemize}
\item $\varphi(w)\leq \varphi(v)$ for $w \in T$ and $v \in T_{w1}$ and
\item $\varphi(w)< \varphi(v)$ for $w \in T$ and $v \in T_{w2}$.
\end{itemize}
We call a map $\varphi: T \to \{1,\ldots , m\}$
a \defit{left-strictly-increasing labeling} if 
\begin{itemize}
\item $\varphi(w)< \varphi(v)$ for $w \in T$ and $v \in T_{w1}$ and
\item $\varphi(w)\leq \varphi(v)$ for $w \in T$ and $v \in T_{w2}$.
\end{itemize}
We call a map $\varphi: T \to \{1,\ldots , m\}$
a \defit{binary-searching labeling} if 
\begin{itemize}
\item $\varphi(w)\geq \varphi(v)$ for $w \in T$ and $v \in T_{w1}$ and
\item $\varphi(w)< \varphi(v)$ for $w \in T$ and $v \in T_{w2}$.
\end{itemize}
\end{definition}

For example, let $T=\{0,1,2,11,21,22\}$. We write a labeling $\varphi$
on $T$ as the diagram
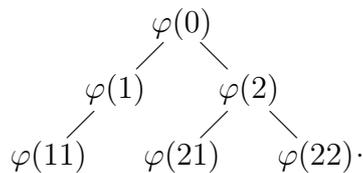
\begin{gather*}
\begin{picture}(130,50)(-65,-50)
\put(0,0){\makebox(0,0)[c]{$\varphi(0)$}}
\put(-25,-25){\makebox(0,0)[c]{$\varphi(1)$}}
\put(25,-25){\makebox(0,0)[c]{$\varphi(2)$}}
\put(-50,-50){\makebox(0,0)[c]{$\varphi(11)$}}
\put(0,-50){\makebox(0,0)[c]{$\varphi(21)$}}
\put(50,-50){\makebox(0,0)[c]{$\varphi(22)$}}
\put(-7,-7){\line(-1,-1){10}}
\put(7,-7){\line(1,-1){10}}
\put(-33,-33){\line(-1,-1){10}}
\put(17,-33){\line(-1,-1){10}}
\put(33,-33){\line(1,-1){10}}
\end{picture}.
\end{gather*}
In this notation, the labelings
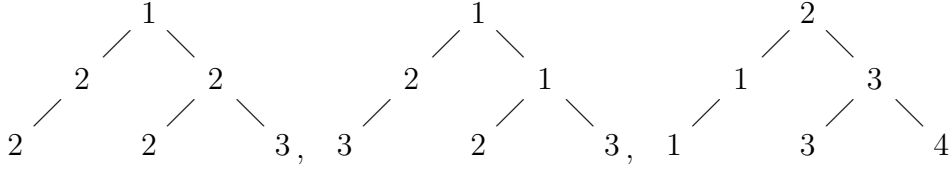
\begin{align*}
\begin{picture}(110,60)(-55,-55)
\put(0,0){\makebox(0,0)[c]{$1$}}
\put(-25,-25){\makebox(0,0)[c]{$2$}}
\put(25,-25){\makebox(0,0)[c]{$2$}}
\put(-50,-50){\makebox(0,0)[c]{$2$}}
\put(0,-50){\makebox(0,0)[c]{$2$}}
\put(50,-50){\makebox(0,0)[c]{$3$}}
\put(-7,-7){\line(-1,-1){10}}
\put(7,-7){\line(1,-1){10}}
\put(-33,-33){\line(-1,-1){10}}
\put(17,-33){\line(-1,-1){10}}
\put(33,-33){\line(1,-1){10}}
\end{picture}&,&
\begin{picture}(110,60)(-55,-55)
\put(0,0){\makebox(0,0)[c]{$1$}}
\put(-25,-25){\makebox(0,0)[c]{$2$}}
\put(25,-25){\makebox(0,0)[c]{$1$}}
\put(-50,-50){\makebox(0,0)[c]{$3$}}
\put(0,-50){\makebox(0,0)[c]{$2$}}
\put(50,-50){\makebox(0,0)[c]{$3$}}
\put(-7,-7){\line(-1,-1){10}}
\put(7,-7){\line(1,-1){10}}
\put(-33,-33){\line(-1,-1){10}}
\put(17,-33){\line(-1,-1){10}}
\put(33,-33){\line(1,-1){10}}
\end{picture}&,&
\begin{picture}(110,60)(-55,-55)
\put(0,0){\makebox(0,0)[c]{$2$}}
\put(-25,-25){\makebox(0,0)[c]{$1$}}
\put(25,-25){\makebox(0,0)[c]{$3$}}
\put(-50,-50){\makebox(0,0)[c]{$1$}}
\put(0,-50){\makebox(0,0)[c]{$3$}}
\put(50,-50){\makebox(0,0)[c]{$4$}}
\put(-7,-7){\line(-1,-1){10}}
\put(7,-7){\line(1,-1){10}}
\put(-33,-33){\line(-1,-1){10}}
\put(17,-33){\line(-1,-1){10}}
\put(33,-33){\line(1,-1){10}}
\end{picture}
\end{align*}
on $T$ are
a right-strictly-increasing labeling,
a left-strictly-increasing labeling
and a binary-searching labeling,  respectively.

The inverse image $\varphi^{-1}(\{1,\ldots,n+1\})$ 
of a right-strictly-increasing labeling $\varphi$
is the tree obtained
from the inverse image $\varphi^{-1}(\{1,\ldots,n\})$
by adding  some nodes right-strictly.
Hence we identify right-strictly-increasing labelings 
with sequences $(\emptyset=T^{0},T^{1},\ldots, T^{m})$ of $m+1$ trees
such that
$T^{i+1}$ is obtained from  $T^{i}$ by adding  some nodes right-strictly
for each $i$.

Similarly, we identify left-strictly-increasing labelings  
with sequences $(\emptyset=T^{0},T^{1},\ldots, T^{m})$ of $m+1$ trees
such that
$T^{i+1}$ is obtained from  $T^{i}$ by adding  some nodes left-strictly
for each $i$.

For a binary-searching labeling $\varphi_{m}:T \to \{1,\ldots, m\}$,
by the definition of binary-searching labeling,
the inverse image $\varphi_{m}^{-1}(\{m\})$
equals $\{w_{T,1},\ldots,w_{T,k}\}$ for some $k$.
We can obtain a binary-searching labeling 
$\varphi_{m-1}: T\ominus \varphi_{m}^{-1}(\{m\}) \to \{1,\ldots, m-1\}$
from $\varphi_{m}$ by evacuating $k$ nodes $\varphi_{m}^{-1}(\{m\})$
together with their labels. 
Hence we identify  binary-searching labelings
with sequences $(\emptyset=T^{0},T^{1},\ldots, T^{m})$ of $m+1$ trees
such that 
$D_{k_i} T^i = T^{i-1}$ for some $k_1$, $k_2$, \dots ,$k_m$.

For a labeling $\varphi$ from $T$ to $\{1,\ldots,m\}$,
we define $t^{\varphi}=\prod_{w\in T} t_{\varphi(w)}$.
For a tree $T$, it follows that
\begin{align*}
\SU[U]{T}{\emptyset}(t_1,\ldots,t_n)
&=\sum_{\varphi} t^{\varphi},\\
\SU[U']{T}{\emptyset}(t_1,\ldots,t_n)
&=\sum_{\phi} t^{\phi},\\
\SD[D]{T}{\emptyset}(t_1,\ldots,t_n)
&=\sum_{\psi} t^{\psi},
\end{align*}
where 
the first sum is over all right-strictly-increasing labelings $\varphi$ on $T$,
the second sum is over all left-strictly-increasing labelings $\phi$ on $T$,
and
the last sum is over all binary-searching labelings $\psi$ on $T$.

\newpage
\thispagestyle{empty}

\newpage
\thispagestyle{empty}


\begin{thebibliography}{99}
\bibitem{BMSW}
Bergeron, Mykytiuk, Sottile and van Willigenburg, 
\emph{Non-commutative Pieri operators on posets},  
J. Combin. Th. Ser. A, 
\textbf{91}, % No. 1/2,
(2000), 84--110.


\bibitem{FomRSK}
S. Fomin,
\emph{Generalized Robinson-Schensted-Knuth correspondence},
Zariski Nauchn. Sem. Leningrad. Otdel. Mat. Inst. Steklov. (LOMI)
\textbf{155}
(1986),
156--175, 195 (Russian);
English transl., 
J. Soviet Math.
41(1988),
979--991.

\bibitem{FomD}
S. Fomin, 
\emph{Duality of graded graphs}, 
J. Algebraic Combin. 
\textbf{3}
(1994),
357--404.

\bibitem{FomS}
S. Fomin,
\emph{Schensted algorithms for dual graded graphs},
J. Algebraic Combin.
\textbf{4}
(1995),
5--45.


\bibitem{FomK}
S. Fomin, 
\emph{Schur operators and Knuth correspondences}, 
J. Combin. Theory
Ser. A 
\textbf{72} 
(1995), 
277--292.

\bibitem{lam}
T. Lam,
\emph{A Combinatorial Generalization of the Boson-Fermion Correspondence},
Math. Res. Lett.
\textbf{13}
(2006),
no.3,
377--329.

\bibitem{tree}
Y. Numata,
\emph{An example of generalized Schur operators involving planar binary
trees},
preprint,
{\tt arXiv:math.CO/0609376}.

\bibitem{Rob}
T. Roby,
\emph{Applications and extensions of Fomin's generalization of the Robinson-Schensted 
 correspondence to differential posets},
Ph.D. thesis,
M.I.T.,
1991.


\bibitem{StaD}
R. Stanley,
\emph{Differential posets},
J. Amer. Math. Soc.,
\textbf{1}
(1988),
919--961.

\bibitem{StaV}
R. Stanley,
\emph{Variations on differential posets},
Invariant theory and tableaux
(Stanton,D.,ed.),
IMA volumes in mathematics
and its applications,
Springer-Verlag,
New York,
145--165.

\end{thebibliography}
\end{document}